\documentclass[letterpaper, preprint, paper,11pt]{IEEEtran}	

\usepackage{bm}
\usepackage{amsmath}
\usepackage{graphicx}
\usepackage{subfigure}
\usepackage{cite}
\usepackage{footnpag}			      	
\usepackage{amsmath} 
\usepackage{amssymb}
\usepackage{amsthm}
\usepackage{bm}
\usepackage{mathrsfs}
\usepackage{hyperref}
\usepackage{float}
\usepackage{accents}

\newcommand{\bma}[1]{\left[\begin{array}{#1}}
\newcommand{\ema}{\end{array}\right]}
\newcommand{\utimes}{ {\raisebox{-0.6ex}{ \kern-1.0ex\raisebox{0.6ex}{ \small $\mathsf{v}$}}} } %

\DeclareMathAlphabet{\mbf}{OT1}{ptm}{b}{n} 

\def\fdotb{{\raisebox{-0.6ex}{ \kern0.2ex\raisebox{0.8ex}{\tiny $\hspace*{-1ex}\circ$}}}}

\def\fddotb{{\raisebox{-0.6ex}{ \kern0.2ex\raisebox{0.8ex}{\tiny $\hspace*{-1ex}\circ\circ$}}}}

\newcommand{\beq}{\begin{equation}}
\newcommand{\eeq}{\end{equation}}
\newcommand{\bdis}{\begin{displaymath}}
\newcommand{\edis}{\end{displaymath}}
\newcommand{\beqarray}{\begin{eqnarray}}
\newcommand{\eeqarray}{\end{eqnarray}}
\newcommand{\beqarraynn}{\begin{eqnarray*}}
\newcommand{\eeqarraynn}{\end{eqnarray*}}

\begin{document}

\title{Attitude Control of a 2U Cubesat by Magnetic and Air Drag Torques}

\author{Richard Sutherland$^1$\thanks{1. PhD Candidate, Department of Aerospace Engineering, University of Michigan, 1320 Beal Ave, Ann Arbor, MI, 48109-2140.},  
Ilya Kolmanovsky$^2$\thanks{2. Professor, Department of Aerospace Engineering, University of Michigan, 1320 Beal Ave, Ann Arbor, MI, 48109-2140.},
and Anouck Girard$^3$\thanks{3. Associate Professor, Department of Aerospace Engineering, University of Michigan, 1320 Beal Ave, Ann Arbor, MI, 48109-2140.}
}

\maketitle{}

\begin{abstract}
This paper describes the development of a magnetic attitude control subsystem for a 2U cubesat. Due to the presence of gravity gradient torques, the satellite dynamics are open-loop unstable near the desired pointing configuration. Nevertheless the linearized time-varying system is completely controllable, under easily verifiable conditions, and the system's disturbance rejection capabilities can be enhanced by adding air drag panels exemplifying a beneficial interplay between hardware design and control. In the paper, conditions for the complete controllability for the case of a magnetically controlled satellite with passive air drag panels are developed, and simulation case studies with the LQR and MPC control designs applied in combination with a nonlinear time-varying input transformation are presented to demonstrate the ability of the closed-loop system to satisfy mission objectives despite disturbance torques.
\end{abstract}

\section{Nomenclature}
\label{sec:nom}
\begin{tabbing}
  XXX \= \kill
  $\mathbf{I}_n$ \> \hspace*{4 mm} n-by-n Identity matrix \\
  $\mathbf{0}_{m \times n}$ \> \hspace*{4 mm} m-by-n Zero matrix \\
  $\mathbf{J}^{\mathcal{B}c}_b$ \> \hspace*{4 mm} Inertia matrix of body $\mathcal{B}$ about center of \\
  \> \hspace*{4 mm} mass, resolved in frame $b$ \\
  $J$ \> \hspace*{4 mm} Optimal control cost functional \\
  $\boldsymbol{\accentset{\rightharpoonup}{a}}$ \> \hspace*{4 mm} Physical vector \\
  $\boldsymbol{a}_f$ \> \hspace*{4 mm} Physical vector resolved in frame $f$ \\
  $\boldsymbol{\hat{a}}$ \> \hspace*{4 mm} Unit vector \\
  $\mathbf{b}$ \> \hspace*{4 mm} Magnetic field vector \\
  $\mathbf{x}$ \> \hspace*{4 mm} State vector \\
  $\boldsymbol{\tau}$ \> \hspace*{4 mm} External torque vector \\
  $\phi$ \> \hspace*{4 mm} Roll angle \\
  $\theta$ \> \hspace*{4 mm} Pitch angle \\
  $\psi$ \> \hspace*{4 mm} Yaw angle \\  
  $\mathbf{O}$ \> \hspace*{4 mm} Orientation matrix \\
  $S[\cdot]$ \> \hspace*{4 mm} Skew-symmetric matrix operator \\
  $c_{(\cdot)}$ \> \hspace*{4 mm} cosine operator \\
  $s_{(\cdot)}$ \> \hspace*{4 mm} sine operator \\
  \textit{Subscripts}\\
  $g$ \> \hspace*{4 mm} Inertial frame \\
  $L$ \> \hspace*{4 mm} Local-Vertical/Local-Horizontal (LVLH) frame \\
  $b$ \> \hspace*{4 mm} Body-fixed principal (BFP) frame \\
\end{tabbing}

\section{Introduction}
\label{sec:intro}
\begin{figure}
    \centering
    \includegraphics[scale = 0.7]{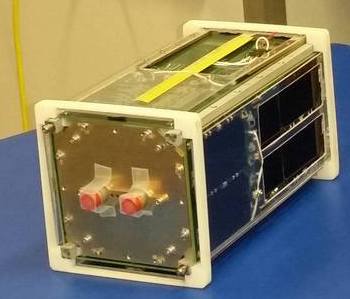}
    \caption{Top face of one of the 2U cubesats, featuring the INMS sensor array. This face must be kept pointed within a 20$^\circ$ half-angle cone of the satellite's orbital velocity vector for the array to function correctly. Two of the air drag panels are visible in their undeployed configurations.}
    \label{fig:cubesat}
\end{figure}

This paper considers the development of an attitude control system, with application to two of the QB50 satellites designed to conduct a survey of the upper atmosphere at low-Earth orbit altitudes\cite{QB50,news}. See Figure \ref{fig:cubesat}. The primary enabler for this survey mission is a constellation of forty 2U cubesats, each equipped with an ion-and-neutral mass spectrometer (INMS) sensor mounted to one of the 1U faces. To function correctly, this sensor must be kept pointed to within a $20^{\circ}$ half-angle cone of the velocity vector. The key challenge in maintaining this attitude is that it corresponds to the cubesat being near a gravity gradient unstable equilibrium. The purpose of this paper is to demonstrate controllability of the linearized time-varying dynamic system and to design a controller for the attitude of a 2U cubesat using first magnetic torque rod actuators alone and in combination with a hardware modification that involves an additional set of four air drag panels. The passive air drag panels are introduced to enhance, in combination with a magnetic rod actuator controller, the satellite's stability and disturbance rejection characteristics.

The cubesat kinematics are expressed through an Euler angle parametrization. The dynamics are characterized using Euler's equation, and incorporate the effects of gravity gradient torque. Controlling the satellite's attitude via magnetic actuators alone is an attractive option; magnetic rods are compact in size, have no moving parts, and consume only electricity, which can be supplied by batteries and solar panels. The drawback to purely magnetic actuation, however, is that the system is instantaneously underactuated due to the inability to exert a torque parallel to the direction of the magnetic field vector. However, unlike underactuated systems involving reaction wheel or thruster failures \cite{NHM1,NHM2,Crouch,Petersen}, in this system the unactuated axis is not fixed with respect to the body but rotates as the satellite traverses its orbit. System controllability is obtained by taking advantage of the time-varying unactuated axis in tandem with the gravity gradient torque.

Much prior work has been done to investigate the use of magnetic torques in spacecraft attitude control. However, that work has largely focused on the use of magnetic torque in spin-stabilized spacecraft \cite{Ferreira,Rodden} or in gravity gradient stabilized spacecraft \cite{Blanke,Steyn}. The satellite considered in this paper is not gravity gradient stabilized; in fact, the desired attitude corresponds to a gravity gradient unstable equilibrium, so a control law is needed that establishes pointing despite being hindered, rather than aided, by the gravity gradient torque. While destablizing in our configuration, it is interesting that the gravity gradient torque, on the other hand, facilitates the satellite controllability and the overcoming of the effects of underactuation. Previous efforts were also made to exploit the quasi-periodicity of the magnetic field in controller design \cite{Psiaki1}. These approaches typically pursue time-averaged solutions of the changing magnetic field to precompute control gains offline; however, these solutions can grow less accurate over time and may require that the satellite be sent updated time-averaged parameters. This work uses magnetic field readings to calculate control gains online. Previous work on passive aerodynamic stability treated small-area drag surfaces \cite{Psiaki2} with a tendency to twist and deform, whereas our satellite is equipped with larger drag plates that should be less prone to performance-degrading deformation. In the conference paper \cite{SFM}, we have presented simulation results for the LQR controller with a design based on a model that did not account for the presence of air drag panels. In this paper, the LQR controller design is based on a linearized, discrete-time model that accounts for the effect of the air drag panels, for which the controllability analysis results are also established. As in \cite{SFM}, before applying the LQR controller and linearizing the model, we use a nonlinear state and control transformation from \cite{Lovera}. The present paper also contains other developments, discussions and details not present in \cite{SFM}. In particular, we present simulation results for the case of an MPC controller that is capable of enforcing control constraints.

This paper is organized as follows. In Section \ref{sec:eom}, equations of motion are introduced. The control analysis and LQR-based controller design for the case of the satellite without air drag panels is introduced in Section \ref{sec:ctrl_law}, with closed-loop simulation results reported in Section \ref{sec:full_sim}. The control analysis and LQR-based control design for the satellite with the added drag panels is presented in Section \ref{sec:darts}. The development of an MPC-based controller that extends the LQR design to be able to handle the control magnitude constraints is given in Section \ref{sec:MPC}. Conclusions are drawn in Section \ref{sec:concl}.

\section{Equations of Motion}
\label{sec:eom}
Our first step in developing a controller is to derive the equations of motion (EOMs) for the system. We choose a body-fixed frame such that the $\boldsymbol{\hat{\imath}}_b$-axis aligns with the INMS sensor array, the $\boldsymbol{\hat{k}}_b$-axis aligns with the satellite's radio antenna, and the $\boldsymbol{\hat{\jmath}}_b$-axis completes the right-hand rule; we also assume this frame to be a principal frame and refer to it as body fixed principal (BFP) frame.

\subsection{Euler Angle Attitude Parametrization}
\label{sec:euler}
To stabilize the satellite to the desired attitude, the controller must account for six states: three independent attitude parameters and three independent angular velocity rates. The satellite's inertial measurement unit gives its angular velocity outputs in terms of Euler angle rates, thus, given that the controller is designed to primarily maintain the satellite near the target orientation, we choose to work with an Euler angle parametrization for the kinematics.

The goal of ram pointing is equivalent to aligning the chosen satellite BFP frame with the non-inertial local-vertical/local-horizontal (LVLH) frame with $\boldsymbol{\hat{\imath}}_L$, $\boldsymbol{\hat{\jmath}}_L$, $\boldsymbol{\hat{k}}_L$ being unit vectors, in which $\boldsymbol{\hat{\imath}}_L$ always points along the orbital track and $\boldsymbol{\hat{k}}_L$ points opposite the orbital radius. The direction cosine matrix (DCM) of the satellite's body fixed frame relative to the LVLH frame is represented using a 3-2-1 Euler angle rotation sequence, such that:
\begin{equation}
    \mathbf{O}_{bL} = \mathbf{O}_{1}(\phi)\mathbf{O}_{2}(\theta)\mathbf{O}_{3}(\psi),
\end{equation}
where
\[
    \mathbf{O}_{1}(\phi) =
    \begin{bmatrix}
        1 & 0 & 0 \\ 0 & \cos{\phi} & \sin{\phi} \\ 0 & -\sin{\phi} & \cos{\phi}
    \end{bmatrix} =
    \begin{bmatrix}
        1 & 0 & 0 \\ 0 & c_{\phi} & s_{\phi} \\ 0 & -s_{\phi} & c_{\phi}
    \end{bmatrix},
\]
\[
    \mathbf{O}_{2}(\theta) =
    \begin{bmatrix}
        \cos{\theta} & 0 & -\sin{\theta} \\ 0 & 1 & 0 \\ \sin{\theta} & 0 & \cos{\theta}
    \end{bmatrix} =
    \begin{bmatrix}
        c_{\theta} & 0 & -s_{\theta} \\ 0 & 1 & 0 \\ s_{\theta} & 0 & c_{\theta}
    \end{bmatrix},
\]
\[
    \mathbf{O}_{3}(\psi) =
    \begin{bmatrix}
        \cos{\psi} & \sin{\psi} & 0 \\ -\sin{\psi} & \cos{\psi} & 0 \\ 0 & 0 & 1
    \end{bmatrix} =
    \begin{bmatrix}
        c_{\psi} & s_{\psi} & 0 \\ -s_{\psi} & c_{\psi} & 0 \\ 0 & 0 & 1
    \end{bmatrix},
\]
where we use shorthand $c$ and $s$ to designate cosine and sine of the argument given in the subscript. When all three angles are at zero, which is the target equilibrium, the $\boldsymbol{\hat{\imath}}_b$-axis aligns with the velocity vector and the $\boldsymbol{\hat{k}}_b$-axis points in the nadir direction.

The drawback of 3-2-1 Euler angles is the kinematic singularity at $\cos \theta = 0$, i.e., $\theta = \pm 90^\circ$, which presents known difficulties for the control design \cite{SDC}, thus we require that the pitch angle be within the $-90^\circ < \theta < 90^\circ$ range before beginning actuation. In the event that the initial attitude violates this constraint, we can sidestep the singularity by rotating about each axis by $\pm 180^\circ$, thereby recasting the attitude to an equivalent set of Euler angle parameters. While in theory the pitch angle could still possibly violate the constraint window some time after actuation begins despite initially satisfying it, in simulations this did not lead to any further actuation problems.

\subsection{Kinematics}
\label{sec:kin}
The angular velocity of the body frame relative to an inertial frame can be decomposed as the sum of the intermediate angular velocity physical vectors:
\begin{equation}
    \boldsymbol{\accentset{\rightharpoonup}{\omega}}^{bg} = \boldsymbol{\accentset{\rightharpoonup}{\omega}}^{bL} + \boldsymbol{\accentset{\rightharpoonup}{\omega}}^{Lg}.
\end{equation}
We resolve this expression in the BFP frame to produce the following expression:
\begin{equation}
    \boldsymbol{\omega}_{b}^{bg} = \boldsymbol{\omega}_{b}^{bL} + \boldsymbol{\omega}_{b}^{Lg} = \boldsymbol{\omega}_{b}^{bL} + \mathbf{O}_{bL} \boldsymbol{\omega}_{L}^{Lg}.
\end{equation}
The LVLH frame rotates, relative to an inertial frame, at a rate given by the orbital motion $n$. This rate is constant due to the circular orbit assumption. In the chosen coordinate system, $\boldsymbol{\omega}_{L}^{Lg} = \left[ 0 \hspace*{1 mm} - \hspace*{-1mm}n \hspace*{2 mm} 0 \right] ^{T}$. Then,
\begin{equation}
    \begin{bmatrix}
        \omega_1 \\ \omega_2 \\ \omega_3
    \end{bmatrix}
    = C_{\phi \theta}
    \begin{bmatrix}
        \dot{\phi} \\ \dot{\theta} \\ \dot{\psi}
    \end{bmatrix}
    + \mathbf{O}_{bL}
    \begin{bmatrix} 0 \\ -n \\ 0 \end{bmatrix},
\label{eqn:ang_vel}
\end{equation}
where
\[
    C_{\phi \theta} =
    \begin{bmatrix}
        1 & 0 & -s_\theta \\
        0 & c_\phi & s_\phi c_\theta \\
        0 & -s_\phi & c_\phi c_\theta
    \end{bmatrix},
\]
and
\[
    \mathbf{O}_{bL} =
    \begin{bmatrix}
        c_\theta c_\psi & c_\theta s_\psi & -s_\theta \\
        s_\phi s_\theta c_\psi - c_\phi s_\psi & s_\phi s_\theta s_\psi + c_\phi c_\psi & s_\phi c_\theta \\
        c_\phi s_\theta c_\psi + s_\phi s_\psi & c_\phi s_\theta s_\psi - s_\phi c_\psi & c_\phi c_\theta
    \end{bmatrix}.
\]
We invert (\ref{eqn:ang_vel}) to solve for the Euler angle rates:
\begin{equation}
    \begin{bmatrix}
        \dot{\phi} \\ \dot{\theta} \\ \dot{\psi} 
    \end{bmatrix}
    = C^{-1}_{\phi \theta}
    \left(
    \begin{bmatrix}
        \omega_1 \\ \omega_2 \\ \omega_3
    \end{bmatrix}
    + n
    \begin{bmatrix}
        c_\theta s_\psi \\
        s_\phi s_\theta s_\psi + c_\phi c_\psi \\
        c_\phi s_\theta s_\psi - s_\phi c_\psi
    \end{bmatrix}
    \right),
\label{eqn:kinematics}
\end{equation}
where
\[
    C^{-1}_{\phi \theta} = \left( \frac{1}{c_\theta} \right)
    \begin{bmatrix}
        c_\theta & s_\phi c_\theta & c_\phi s_\theta \\
        0 & c_\phi c_\theta & -s_\phi c_\theta \\
        0 & s_\phi & c_\phi
    \end{bmatrix}.
\]

\subsection{Dynamics}
\label{sec:dyn}
Having derived the equations for the kinematics, we turn to the dynamics, which can be modeled with Euler's equation,
\begin{equation}
    \mathbf{J}^{\mathcal{B}c}_b\dot{\boldsymbol{\omega}}^{bg}_b + S[\boldsymbol{\omega}^{bg}_b]\mathbf{J}^{\mathcal{B}c}_b\boldsymbol{\omega}^{bg}_b = \boldsymbol{\tau}^{\mathcal{B}c}_b,
\label{eqn:dynamics}
\end{equation}
where all quantities have been resolved in BFP frame, $\boldsymbol{\tau}^{\mathcal{B}c}_b$ represents the external torque on body $\mathcal{B}$ about its center of mass, and the matrix $S[\boldsymbol{\omega}^{bg}_b]$ denotes the skew-symmetric matrix formed of the components of $\boldsymbol{\omega}^{bg}_b$ and given by
\begin{equation}
    S[\boldsymbol{\omega}^{bg}_b] = S\left(
    \begin{bmatrix}
        \omega_1 \\ \omega_2 \\ \omega_3
    \end{bmatrix} \right) =
    \begin{bmatrix}
        0 & -\omega_3 & \omega_2 \\ \omega_3 & 0 & -\omega_1 \\ -\omega_2 & \omega_1 & 0
    \end{bmatrix}.
\end{equation}
\noindent
Since the BFP frame is a principal frame, $\mathbf{J}^{\mathcal{B}c}_b$ is
\begin{equation}
    \mathbf{J}^{\mathcal{B}c}_b =
    \begin{bmatrix}
        J_1 & 0 & 0 \\
        0 & J_2 & 0 \\
        0 & 0 & J_3
    \end{bmatrix},
\label{eqn:J_BFP}
\end{equation}
where $J_1, J_2,$ and $J_3$ are the principal moments of inertia. The torque acting on the body can be further decomposed into the magnetic control torque, gravity gradient torque, and disturbance torque,
\begin{equation}
\boldsymbol{\tau}^{\mathcal{B}c}_b = \boldsymbol{\tau}^{\mathcal{B}c,mt}_{b} + \boldsymbol{\tau}^{\mathcal{B}c,gg}_{b} + \boldsymbol{\tau}^{\mathcal{B}c,dist}_{b}.
\end{equation}

\subsubsection{Magnetic Torque}
\label{sec:mt}
The net magnetic dipole generated by the torque rods interacts with the Earth's magnetic field to produce a torque according to the following law:
\begin{equation}
    \boldsymbol{\accentset{\rightharpoonup}\tau}^{\mathcal{B}c,mt} = \boldsymbol{\accentset{\rightharpoonup}m} \times \boldsymbol{\accentset{\rightharpoonup}b},
\label{eqn:mt}
\end{equation}
where $\boldsymbol{\accentset{\rightharpoonup}b}$ denotes the external magnetic field vector and $\boldsymbol{\accentset{\rightharpoonup}m}$ denotes the magnetic moment. Resolving the torque in the BFP frame, we obtain
\begin{equation}
    \boldsymbol{\tau}_{b}^{\mathcal{B}c,mt} = -S[\boldsymbol{b}_{b}]\mathbf{m}_b,
\end{equation}
where we have replaced the vector cross-product operation with the equivalent skew-symmetric matrix representation,
\begin{equation}
    S\left(
    \begin{bmatrix}
        b_1 \\ b_2 \\ b_3
    \end{bmatrix} \right) =
    \begin{bmatrix}
        0 & -b_3 & b_2 \\ b_3 & 0 & -b_1 \\ -b_2 & b_1 & 0
    \end{bmatrix}.
\label{eqn:skew}
\end{equation}
Note that rank$(S[\boldsymbol{b}_b]) = 2$, which affirms that the system is instantaneously underactuated with magnetic torque alone.

\subsubsection{Gravity Gradient Torque}
\label{sec:gg}
The Earth's gravitational field exerts a force on the satellite that can be modeled by an inverse-square distance law. Hence, the force is slightly greater on the portions of the satellite that are closer to the Earth than on those portions that are further away. This differential, though small, produces a torque on the satellite's body given by \cite{Wie}
\begin{equation}
    \boldsymbol{\tau}_{b}^{\mathcal{B}c,gg} = 3 n^{2} S[\boldsymbol{\hat{r}}_{b}] J_{b}^{\mathcal{B}c} \boldsymbol{\hat{r}}_{b},
\end{equation}
where a circular orbit is assumed for the satellite. The radial unit vector $\boldsymbol{\hat{r}}_{b}$ points opposite to the $\boldsymbol{\hat{k}}_L$ vector; hence, $\boldsymbol{\hat{r}}_b$ can be expressed as
\[
    \boldsymbol{\hat{r}}_b = \mathbf{O}_{bL} (-\boldsymbol{\hat{k}}_{L}) = \mathbf{O}_{bL}
    \begin{bmatrix} 0 \\ 0 \\ -1 \end{bmatrix}.
\]
Thus,
\begin{equation}
    \boldsymbol{\tau}_{b}^{\mathcal{B}c,gg} = 3 n^{2}
    \begin{bmatrix}
        \left(J_2 - J_3\right) c_{\phi} s_{\phi} c_{\theta}^{2} \\
        \left(J_3 - J_1\right) c_{\phi} c_{\theta} s_{\theta} \\
        \left(J_1 - J_2\right) s_{\phi} c_{\theta} s_{\theta}
    \end{bmatrix}.
\end{equation}
There exist multiple configurations in which this torque vanishes, two stable $\left( c_\theta = 0 \right)$ and several unstable $(s_{\theta} = s_{2\phi} = 0)$. As the desired equilibrium attitude for the satellite is close to one of the unstable configurations, the gravity gradient torque tends to destabilize the satellite's attitude, and thus becomes a significant factor in our treatment of the system dynamics.

\section{Control Analysis and Control Law Design for the Satellite Without Drag Panels}
\label{sec:ctrl_law}
The magnetic torque cross-product expression (\ref{eqn:mt}) indicates that the component of the dipole parallel to the magnetic field generates zero torque. Thus, to obtain maximum control torque, we seek a dipole moment that minimizes the magnitude of the projection $\boldsymbol{\accentset{\rightharpoonup}m} \cdot \boldsymbol{\accentset{\rightharpoonup}b} = \boldsymbol{m}_b^T \boldsymbol{b}_b$. Following Lovera and Astolfi \cite{Lovera}, we prescribe a dipole moment of the form:
\begin{equation}
    \mathbf{m}_b = - \left( \frac{S[\boldsymbol{b}_b]}{\boldsymbol{b}_b^T \boldsymbol{b}_b} \right) \mathbf{u},
\end{equation}
where $\mathbf{u} \in \mathbb{R}^3$ is a control input vector. This dipole law constrains $\mathbf{m}_b$ to be perpendicular to $\boldsymbol{b}_b$, thus fixing $\mathbf{m}_b^T \boldsymbol{b}_b = 0$. The magnetic input torque can then be expressed as
\begin{equation}
    \boldsymbol{\tau}_b^{\mathcal{B}c,mt} = \left( \frac{S[\boldsymbol{b}_b]S[\boldsymbol{b}_b]}{\boldsymbol{b}_b^T \boldsymbol{b}_b} \right) \mathbf{u}.
\end{equation}
Lovera and Astolfi \cite{Lovera} use a PD control law to prescribe the new control input vector $\mathbf{u}$; we choose instead to apply Linear Quadratic Regulator (LQR) theory to obtain the controller as it can be applied systematically to different spacecraft configurations.

\subsection{Linearized Equations of Motion}
\label{sec:lin_eom}
We first linearize the EOMs about the desired equilibrium state $\mathbf{x}_{eq} = [0, 0, 0, 0, -n, 0]^T$ to get linearized EOMs in the form
\[
    \mathbf{\dot{x}} = A_c \mathbf{x} + B_c(t) \mathbf{u},
\]
with
\begin{equation}
    A_{c} =
    \begingroup 
    \setlength\arraycolsep{2.5pt}
    \begin{bmatrix}
        0 & 0 & n & 1 & 0 & 0 \\
        0 & 0 & 0 & 0 & 1 & 0 \\
        -n & 0 & 0 & 0 & 0 & 1 \\
        -3n^{2} J_{23} & 0 & 0 & 0 & 0 & -n J_{23} \\
        0 & 3 n^{2} J_{31} & 0 & 0 & 0 & 0 \\
        0 & 0 & 0 & -n J_{12} & 0 & 0
    \end{bmatrix}
    \endgroup ,
\label{eqn:full_Ac}
\end{equation}
where $J_{12} := (J_1 - J_2)/J_3$, $J_{23} := (J_2 - J_3)/J_1$, and $J_{31} := (J_3 - J_1)/J_2$, and
\begin{equation}
    B_{c}(t) = 
    \begin{bmatrix}
        \mathbf{0}_{3 \times 3} \\ \\
        \left(\mathbf{J}_b^{\mathcal{B}c}\right)^{-1} \frac{S[\boldsymbol{b}_b(t)]S[\boldsymbol{b}_b(t)]}{\boldsymbol{b}_b^{T}(t) \boldsymbol{b}_b(t)}
    \end{bmatrix}.
\label{eqn:full_Bc}
\end{equation}
The complete derivation of the linearized EOMs is found in Appendix A. Note that $\mathbf{b}_b(t)$ in (\ref{eqn:full_Bc}) is ideally the nominal magnetic field at the linearization point. However, in our implementation, we use the measured magnetic field values in (\ref{eqn:full_Bc}), which are determined by the onboard magnetomer readings. While this makes little difference in terms of model accuracy near the linearization point, this approach allows us to implement the controller with gains computed online without strong coupling to the nominal orbital position or the need to either store nominal values of $\mathbf{b}_b$ or compute them offline and store them in ROM. This also improves the robustness in case of deployment errors or orbit decaying due to the influence of the air drag. In either case, $\mathbf{b}_b$ in (\ref{eqn:full_Bc}) depends on time $t$, as does $B_c(t)$. Note that the matrices $A_c, B_c(t)$ at any fixed time instant $t$ do not constitute a controllable pair, e.g., they violate the controllability rank condition for time-invariant systems.

\subsubsection{Controllability of the Time-Varying Linearized System}
\label{sec:lin_con}
Having derived the linearized EOMs, we now demonstrate that the linearized system, which is time-varying, is controllable on any time interval of non-zero length. Yang \cite{Yang} reduces the problem of complete controllability by magnetic torque rods to a small number of sufficient conditions:
\begin{itemize}
    \item The satellite is not located on the magnetic Equator.
    \item Assuming the above holds, then
    \begin{enumerate}
        \item $J_2 \neq J_3$,
        \item $6 J_3 \left( J_3 - J_1 \right) \neq J_2 \left( J_1 - J_2 + J_3 \right)$.
    \end{enumerate}
\end{itemize}
The considered cubesat has an inertia matrix of $\mathbf{J}_b^{\mathcal{B}c} = $ diag(3654338, 9060235, 8813148) g$\cdot$mm$^2$. As it is to be ejected from the ISS, it will operate with an initial altitude of 415 km and at an inclination of $51.6^\circ$, with an orbital period of 5570 s. At this inclination, the satellite is away from the plane of the magnetic Equator, and also satisfies the controllability constraints on the $J_i$'s above, confirming that the linearized time-varying system is controllable. This implies that, in the absence of control constraints, there exists a control input that drives the state to the origin over any specified time interval.

\subsection{LQR Control Law Design}
\label{sec:lqr_des}
The control design is based on LQR theory applied to a discrete-time model that is obtained by converting (\ref{eqn:full_Ac})-(\ref{eqn:full_Bc}) to discrete-time. We apply a Zero-Order Hold method to carry out the discretization. Let $t \in \mathbb{R}_{\geq0}$ be the current time instant and $A_c$, $B_c(t)$ defined in (\ref{eqn:full_Ac}), (\ref{eqn:full_Bc}). For $\Delta t > 0$, the discrete-time model predicts the state $\mathbf{x}_{k}$ at time $t+k\Delta t, k \in \mathbb{Z}_{\geq 0}$, according to the following model with the ``frozen-in-time'' magnetic field:
\begin{equation}
\begin{split}
    \mathbf{x}_{k+1} & = A_d \mathbf{x}_{k} + B_d(t) \mathbf{u}_{k}, \\
    A_d & = e^{A_c \Delta t}, \\
    B_d(t) & = - A_c^{-1} \left( \mathbf{I}_6 - A_d \right) B_c(t), \\
    \mathbf{x}_0 &= \mathbf{x}(t).
\end{split}
\label{eqn:zoh}
\end{equation}
The pair $(A_d, B_d(t))$ can be verified to be controllable for all $t$ for our orbit and choices of $\Delta t$ we have used. For the difference equation (\ref{eqn:zoh}), we define the infinite-horizon cost functional $J$:
\begin{equation}
    J(t) = \sum_{k = 0}^{\infty} \left( \mathbf{x}_k^T Q \mathbf{x}_k + \mathbf{u}_k^T R \mathbf{u}_k \right),
\end{equation}
where $R=R^T \in \mathbb{R}^3$ is a positive definite matrix, and $Q=Q^T \in \mathbb{R}^6$ is a positive semi-definite matrix satisfying the usual detectability assumption. Then, the optimal feedback control sequence $\mathbf{u}_k = -K(t) \mathbf{x}_k$ that minimizes $J(t)$ has the solution:
\begin{equation}
    K(t) = \left( R + B_d(t)^T P(t) B_d(t) \right)^{-1} B_d(t)^T P(t) A_d,
\label{eqn:kd}
\end{equation}
\begin{equation}
    P(t) = A_d^T P(t) A_d + Q - A_d^T P(t) B_d(t) K(t).
\label{eqn:dare}
\end{equation}
Note that (\ref{eqn:dare}) can have multiple solutions; the $P(t)$ of interest to us is the unique positive definite solution. Also note that $B_d(t)$ changes throughout the orbit, thus (\ref{eqn:dare}) is solved at different instants $t$ in time and the gain $K(t)$ in (\ref{eqn:kd}) is time-varying. The control $\mathbf{u}(t+\sigma) = K(t) \mathbf{x}(t)$ is applied for $0 \leq \sigma < \Delta t$ and then recomputed. A fast update scheme for the solution of the Algebraic Riccati Equation in response to changes in the magnetic field vector can be defined, see Appendix B.

Though we solve for the optimal control without placing any limits on the solution, in practice there do exist control constraints in the form of magnetic torque rod saturation. We indirectly take these constraints into account in the LQR control formulation; if at least one component violates the constraint, the control input is rescaled such that its largest component is equal to the saturation limit while remaining parallel to the calculated dipole vector. In this way, the dipole vector remains perpendicular to the magnetic field vector, for maximum torque generation. In Section \ref{sec:MPC}, we extend the LQR controller to a Model Predictive Control (MPC) based controller that has the additional capability of being able to explicitly handle these control saturation constraints.

\section{Simulation Results}
\label{sec:full_sim}
After several tuning experiments with the goal of obtaining good performance, the weighting matrices $Q$ and $R$ for the LQR cost functional have been chosen as $Q =$ diag($10^{-8}$, $10^{-8}$, $10^{-8}$, $10^{-4}$, $10^{-4}$, $10^{-4}$) and $R =$ diag($10^8$, $10^8$, $10^8$). The initial attitude of the satellite is $\phi(0) = -35^\circ$, $\theta(0) = -75^\circ$, and $\psi(0) = 75^\circ$. It is estimated that the magnitudes of the post-ejection tumble rates would be approximately $10^\circ$ per second in each axis, thus the satellite is given an initial angular velocity of $\omega_{1}(0) = -10^\circ$/s, $\omega_{2}(0) = 10^\circ$/s, and $\omega_{3}(0) = -10^\circ$/s. The Earth's magnetic field is modeled using the tilted dipole approximation of Wertz \cite{Wertz}, updated with 2015 IGRF-12 \cite{IGRF} coefficients. The controller saturation limit is $u_{max} = 0.1$ A$\cdot$m$^2$. The control goal is to drive the pointing angle, i.e., the angle between $\boldsymbol{\hat{\imath}}_b$ and the satellite's velocity vector, to within the $\pm20^\circ$ constraint.
\begin{figure}
    \centering
    \includegraphics[scale=0.33]{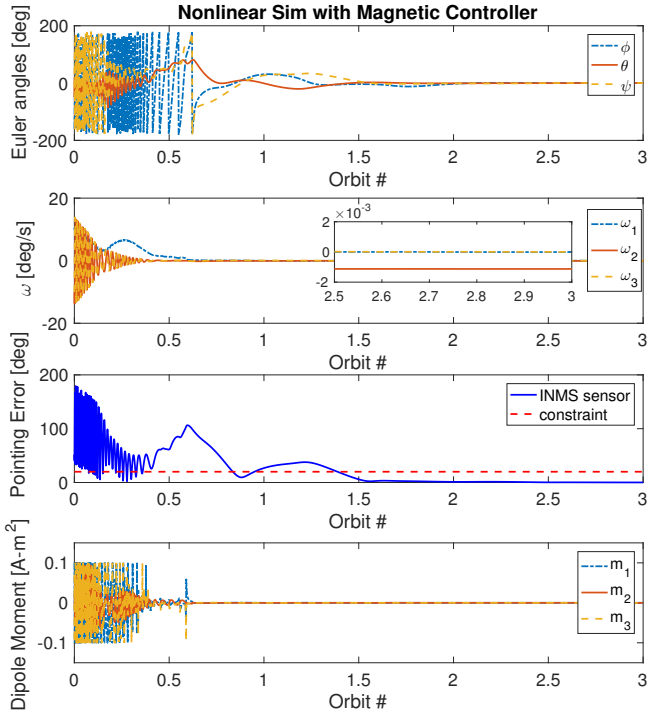}
    \caption{Simulation of the closed-loop system with the LQR control law for the magnetic torque rods with no sensor noise or unmodeled disturbance torques.}
    \label{fig:sim_1}
\end{figure}
The controller successfully achieves the commanded equilibrium, without unwinding. As shown in Figure \ref{fig:sim_1}, the cubesat experiences many rotations about its roll and yaw axes while detumbling, but does reach and remain within the required pointing angle constraint within two and a half orbits and with no unwinding. The pointing error of the INMS sensor approaches zero, indicating that it is correctly oriented.

\section{Secondary Actuation}
\label{sec:darts}
\begin{figure}
    \centering
    \includegraphics[scale=0.32]{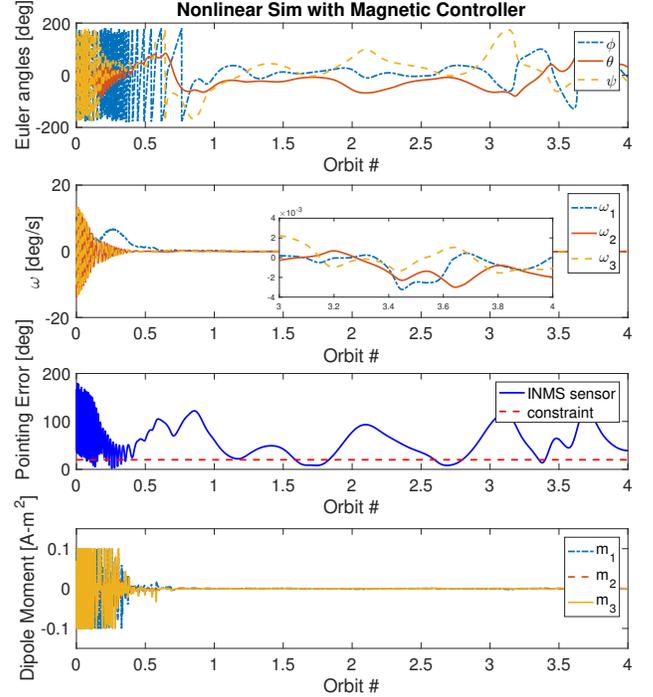}
    \caption{Simulation of the closed-loop system with the LQR control law for the magnetic torque rods with an unmodeled disturbance torque of constant $10^{-8}$ N$\cdot$m magnitude.}
    \label{fig:dist}
\end{figure}

The magnetic torque rod controller works well in the ideal case, but can struggle to maintain the pointing constraint in the presence of unmodeled disturbance torques, as seen in Figure \ref{fig:dist} when a constant magnitude unmodeled disturbance torque is added to the system.
As magnetic rods alone do not appear to provide strong disturbance rejection capability, a solution that takes advantage of hardware and control interplay has been adopted. Specifically, the cubesat was augmented with a set of four drag panels, to supplement the magnetic torque with passive aerodynamic stabilizing torque. We note that in a different application to a flying wind turbine \cite{IK}, a solution that exploits passive aerodynamic stabilization to enhance an underactuated system has also been proposed.
\begin{figure}
    \centering
    \includegraphics[scale=0.52]{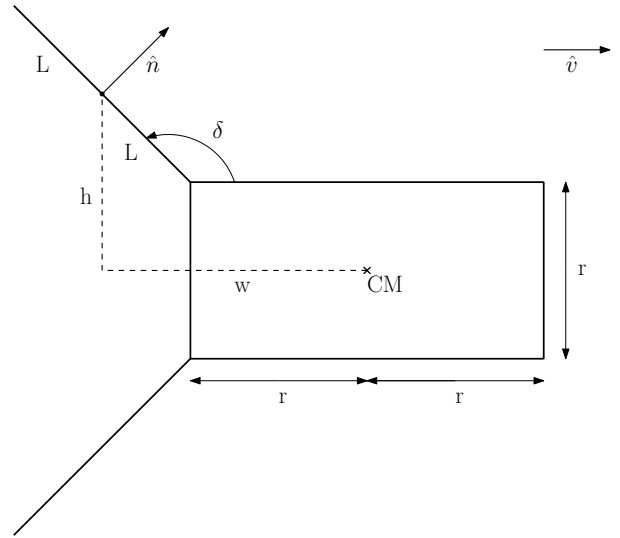}
    \caption{Drag panel model following deployment. The satellite has four such panels but, for clarity, only two panels are depicted here.}
    \label{fig:darts}
\end{figure}

\subsection{Linearization of the Augmented Satellite Dynamics}
\label{sec:panel_lin}
As the panels add additional torque inputs, they change the satellite's dynamics. Thus, in order to apply the LQR controller, the equations of motion must be relinearized to properly reflect these changes. From the circular orbit assumption, the orbital velocity vector satisfies $\boldsymbol{\hat{v}}_L = [ 1, \ 0, \ 0]^T$, and we can apply the orientation matrix $\mathbf{O}_{bL}$ to resolve it in the BFP frame,
\begin{equation}
\begin{split}
    \boldsymbol{\hat{v}}_b &= \begin{bmatrix} c_\theta c_\psi & \sim & \sim \\ s_\phi s_\theta c_\psi - c_\phi s_\psi & \sim & \sim \\ c_\phi s_\theta c_\psi + s_\phi s_\psi & \sim & \sim \end{bmatrix} \begin{bmatrix} 1 \\ 0 \\ 0 \end{bmatrix} \\ &= \begin{bmatrix} c_\theta c_\psi \\ s_\phi s_\theta c_\psi - c_\phi s_\psi \\ c_\phi s_\theta c_\psi + s_\phi s_\psi \end{bmatrix} \approx \begin{bmatrix} 1 \\ \phi \theta - \psi \\ \theta + \phi \psi \end{bmatrix} \approx \begin{bmatrix} 1 \\ -\psi \\ \theta \end{bmatrix},
    \label{eqn:vb}
\end{split}
\end{equation}
where $\boldsymbol{\hat{v}}_b$ is simplified by first a small angle approximation and then the discarding of higher-order terms.

\begin{figure}
    \centering
    \includegraphics[scale=0.28]{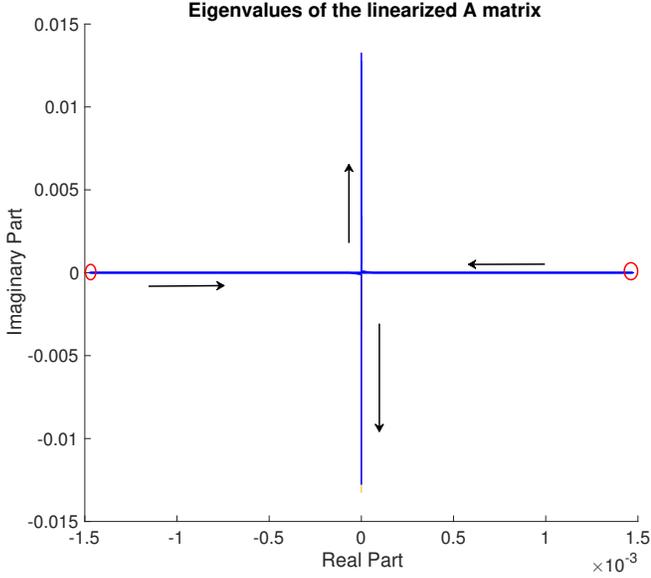}
    \caption{Eigenvalues of the $A$ matrix as a function of panel area. The red circles denote the pair of real eigenvalues that correspond to zero panel area. As panel area increases, this pair of real eigenvalues migrate to the $\jmath \omega$-axis, while the remaining imaginary eigenvalues move outward along this axis.}
    \label{fig:eig}
\end{figure}

\begin{figure}
    \centering
    \includegraphics[scale = 0.3]{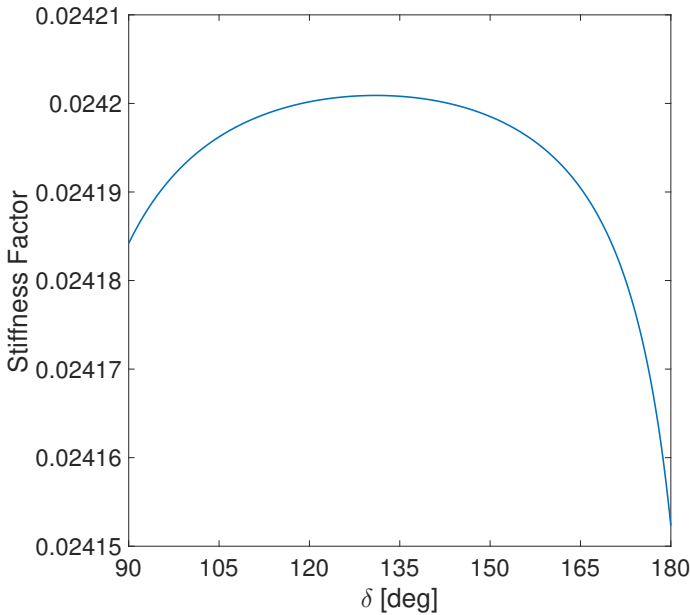}
    \caption{Square root of the magnitude of the smallest-magnitude eigenvalue of the linearized panel system as a function of the panel angle. As all eigenvalues for this particular panel size lie on the $\jmath \omega$-axis, this can be treated as a measure of the ``stiffness'' of the system. For our particular satellite, the panels have their strongest effect at $\delta = 131^\circ$.}
    \label{fig:stiffness}
\end{figure}

Each panel is assumed to contribute a torque of the form
\begin{equation}
    \left(\boldsymbol{\tau}^{ad}_{b}\right) _i = \left(\mathbf{r}_b\right) _i \times \left(\mathbf{F}^{ad}_b\right) _i,
\end{equation}
where $\left(\mathbf{F}^{ad}_b\right) _i$ is the air drag force on the $i^{th}$ panel, $\left(\mathbf{r}_b\right) _i$ is the distance from the satellite's center of mass to the center of the $i^{th}$ panel, and
\begin{equation}
    \boldsymbol{\tau}^{ad}_b = \sum^4_{i=1} \left(\boldsymbol{\tau}^{ad}_{b}\right) _i.
\end{equation}
The force exerted by the air drag acts opposite to the orbital velocity vector, and is modeled as
\begin{equation}
    \begin{split}
        \left(\mathbf{F}^{ad}_b\right) _i &= (0.5 \rho v^2 \mathcal{A}_i C_D) (-\boldsymbol{\hat{v}}_b) \\ &= (0.5 \rho n^2 a^2 \mathcal{A}_i C_D) (-\boldsymbol{\hat{v}}_b),
    \end{split}
\end{equation}
where the circular orbit assumption is used to conclude that $v = na$, with $a$ being the orbital radius, and the effective panel area $\mathcal{A}_i = (\boldsymbol{\hat{n}}_i \cdot \boldsymbol{\hat{v}}) \mathcal{A}$ is equal to the actual panel area scaled by the dot product of the outward facing unit normal vector and the unit velocity vector. Thus, the torque has the following components,
\begin{equation}
    \begin{split}
        \left(\mathbf{r}_b\right)_1 &= \begin{bmatrix} -w \\ 0 \\ -h \end{bmatrix}, \ \
        \left(\mathbf{F}^{ad}_b\right)_1 = n^2 f \left( s_\delta + c_\delta \theta \right) \left( \begin{bmatrix} -1 \\ \psi \\ -\theta \end{bmatrix} \right), \\
        \left(\mathbf{r}_b\right)_2 &= \begin{bmatrix} -w \\ 0 \\ h \end{bmatrix}, \ \
        \left(\mathbf{F}^{ad}_b\right)_2 = n^2 f \left( s_\delta - c_\delta \theta \right) \left( \begin{bmatrix} -1 \\ \psi \\ -\theta \end{bmatrix} \right), \\
        \left(\mathbf{r}_b\right)_3 &= \begin{bmatrix} -w \\ -h \\ 0 \end{bmatrix}, \ \
        \left(\mathbf{F}^{ad}_b\right)_3 = n^2 f \left( s_\delta - c_\delta \psi \right) \left( \begin{bmatrix} -1 \\ \psi \\ -\theta \end{bmatrix} \right), \\
        \left(\mathbf{r}_b\right)_4 &= \begin{bmatrix} -w \\ h \\ 0 \end{bmatrix}, \ \
        \left(\mathbf{F}^{ad}_b\right)_4 = n^2 f \left( s_\delta + c_\delta \psi \right) \left( \begin{bmatrix} -1 \\ \psi \\ -\theta \end{bmatrix} \right),
    \end{split}
\end{equation}
where $f = 0.5 \rho a^2 \mathcal{A} C_D$ and $\delta$ is constrained to be in the interval [90$^\circ$, 180$^\circ$]. Taking the cross products and summing to get an approximation of the air drag torque $\boldsymbol{\tau}^{ad}_b$,
\begin{equation}
    \boldsymbol{\tau}^{ad}_b = 0.5 \rho n^2 a^2 \mathcal{A} C_D \begin{bmatrix} 0 \\ (2 h c_\delta - 4 w s_\delta) \theta \\ (2 h c_\delta - 4 w s_\delta) \psi \end{bmatrix}.
\end{equation}
For the 2U cubesat depicted in Figure \ref{fig:darts}, $w = r - L c_\delta$ and $h = 0.5r + L s_\delta$, thus,
\begin{equation}
\begin{split}
    \boldsymbol{\tau}^{ad}_b &= n^2 f 
    \begin{bmatrix} 
        0 \\ (2 h c_\delta - 4 w s_\delta) \theta \\ (2 h c_\delta - 4 w s_\delta) \psi
    \end{bmatrix} \\ &=
    4 n^2 r f \left(c_\delta - 4 s_\delta + 3 (\frac{L}{r}) s_{2\delta} \right) \begin{bmatrix} 0 \\ \theta \\ \psi \end{bmatrix} = n^2 \Gamma \begin{bmatrix} 0 \\ \theta \\ \psi \end{bmatrix},
\end{split}
\end{equation}
where $\Gamma = 4 r f \left(c_\delta - 4 s_\delta + 3 (\frac{L}{r}) s_{2\delta} \right)$. The linearized dynamic contribution of $\left(\mathbf{J}^{\mathcal{B}c}_b\right)^{-1} \boldsymbol{\tau}^{ad}_b$ then takes the form
\begin{equation}
    \left(\mathbf{J}^{\mathcal{B}c}_b\right)^{-1} d\boldsymbol{\tau}^{ad}_b = 
    \begin{bmatrix}
        0 & 0 & 0 & 0 & 0 & 0 \\
        0 & \frac{n^2 \Gamma}{J_2} & 0 & 0 & 0 & 0 \\
        0 & 0 & \frac{n^2 \Gamma}{J_3} & 0 & 0 & 0
    \end{bmatrix}
    \begin{bmatrix}
        d\phi \\ d\theta \\ d\psi \\ d\omega_1 \\ d\omega_2 \\ d\omega_3
    \end{bmatrix} ,
    \label{eqn:Gamma}
\end{equation}
where $d(\cdot)$ denotes the deviation from the nominal values. Incorporating this contribution into the previously linearized equations of motion $\dot{x} = A x + B u$, $A$ now takes the form
\begin{equation} A = 
    \begingroup 
    \setlength\arraycolsep{2.5pt}
    \begin{bmatrix}
        0 & 0 & n & 1 & 0 & 0 \\
        0 & 0 & 0 & 0 & 1 & 0 \\
        -n & 0 & 0 & 0 & 0 & 1 \\
        3n^2 J_{32} & 0 & 0 & 0 & 0 & n J_{32} \\
        0 & 3n^2 J_{31} + \frac{n^2 \Gamma}{J_2} & 0 & 0 & 0 & 0 \\
        0 & 0 & \frac{n^2 \Gamma}{J_3} & n J_{21} & 0 & 0
    \end{bmatrix}
    \endgroup .
\end{equation}
Notable about this modified $A$ matrix is that, for sufficiently large values of $\Gamma$, all eigenvalues of $A$ lie on the $\jmath \omega$-axis, whereas for $\Gamma = 0$ there exists an unstable positive real eigenvalue, as seen in Figure \ref{fig:eig}. Further, the panel deployment angle directly influences the eigenvalues of the system and its ``stiffness'', i.e., the ability to resist to disturbances. See Figure \ref{fig:stiffness}.

\subsection{LTV Controllability}
\label{sec:ctrb}

The controller uses magnetic rods, thus the $B(t)$ matrix depends on the Earth's magnetic field and the system is time-varying. While for implementation, we rely on measured magnetic field values, for the controllabilty analysis here, we approximate the magnetic field by assuming a tilted dipole model, with periodicity equal to the satellite's orbital period $T$; following Psiaki\cite{Psiaki1}, the magnetic field approximation takes the following form,
\begin{equation}
    \begin{bmatrix}
        b_1(t) \\ b_2(t) \\ b_3(t)
    \end{bmatrix} = \frac{\mu_f}{a^3} 
    \begin{bmatrix}
        \cos(nt) \sin(i_m) \\ -\cos(i_m) \\ 2 \sin(nt) \sin(i_m)
    \end{bmatrix},
\end{equation}
where $\mu_f$ is the strength of the dipole field, $a$ is the semimajor axis, $n$ is the mean motion, $i_m$ is the inclination of the orbit with respect to the magnetic Equator, and $t \in [0, \ T]$ is measured from the crossing of the ascending node of the magnetic Equator. We then show that the LTV system is controllable on the interval $[0, \ T]$, under a few conditions.
\\\\
\textbf{Theorem 1:} The linearized system is controllable on the interval $[0, \ T]$ if the following conditions hold:

1) The satellite's orbital plane is not aligned with the magnetic Equator,

2) $J_2 \neq J_3$,

3) $J_3 \left( 6(J_3 - J_1) + 2 \Gamma \right) \neq J_2 \left(J_1 - J_2 + J_3 - 2 \Gamma \right)$.
\\\\
The proof is similar to the method in Yang \cite{Yang} and is found in Appendix C. Note that if $\Gamma \rightarrow 0$, for example, if the panel area or the atmospheric density were 0, then the last condition specializes to the controllability result for the spacecraft without panels.

\begin{figure}
    \centering
    \includegraphics[scale=0.32]{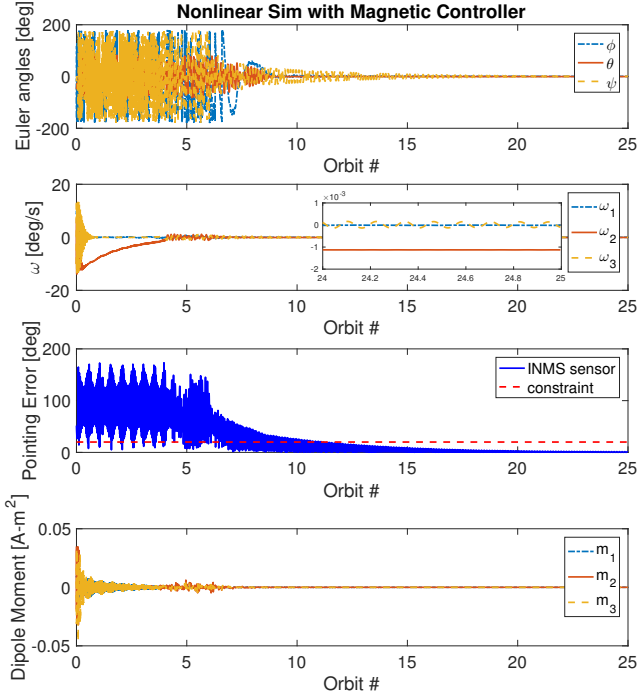}
    \caption{Simulation results using the control design based on the relinearized dynamics, to include the panels, with the same constant disturbance torque as in Figure \ref{fig:dist}. }
    \label{fig:sim_2}
\end{figure}

In Figure \ref{fig:sim_2}, the air drag panels are added to the simulation in Figure \ref{fig:dist}; the convergence is slowed, but the cubesat detumbles properly as the augmented controller does reject the unmodeled disturbance, demonstrating that the additional restorative torques from the panels help to overcome the unmodeled disturbance that destabilized the system previously.

\section{Model Predictive Control}
\label{sec:MPC}

We compare the results of the simulated LQR controller to those of a simulated MPC controller with the same weights as the LQR controller. Unlike the LQR controller, the MPC-based controller has the additional capability of explicitly handling control constraints, such as the magnetic torque rod saturation described in Section \ref{sec:lqr_des}.

\subsection{Discrete Time Conversion}
\label{sec:dtc}
To develop our predictive controller, a discrete-time approximation to the continuous-time dynamics is implemented. The discretization is performed under the assumption that the magnetic field is constant during each control actuation step; this assumption is reasonable, especially over the short prediction horizon, on the order of seconds, that we consider, as the field, while time-varying, is only slowly-varying, with a period of 24 hours.

A zero-order hold, identical to that used in (\ref{eqn:zoh}), is applied to discretize the continuous-time dynamics of the satellite with drag panels system and predict the future state $\mathbf{x}_{k+1}$ according to the ``frozen-in-time'' magnetic field $B_d(t)$.

At each sampling step, the controller then minimizes the now finite-horizon cost functional
\begin{equation}
    J(t) = \mathbf{x}_N^T P(t) \mathbf{x}_N + \sum_{k = 0}^{N-1} \left( \mathbf{x}_k^T Q \mathbf{x}_k + \mathbf{u}_k^T R \mathbf{u}_k \right),
\end{equation}
with prediction horizon $N$ and subject to the discrete-time dynamics in (\ref{eqn:zoh}), as well as to the constraint
\begin{equation}
    |\mathbf{u}_k|_\infty \leq u_{max},
\end{equation}
where $P(t)$ is the unique positive definite solution to the associated Discrete-Time Algebraic Riccati Equation,
\begin{equation}
\begin{split}
    P(t) &= A_d^T P(t) A_d + Q - A_d^T P(t) B_d(t) K(t), \\
    K(t) &= \left( R + B_d(t)^T P(t) B_d(t) \right)^{-1} B_d(t)^T P(t) A_d.
\label{eqn:dare_mpc}
\end{split}
\end{equation}
The controller implements the first control action and then recomputes a new minimizing control sequence at the next sampling time instant.

\subsection{Simulation Results}
\label{sec:sim_mpc}
The MPC controller uses the same weights as in the LQR controller. The discrete-time steps are of length $\Delta t = 4$ sec, and the prediction horizon is held at $N = 5$. All other parameters are identical to those used to generate the simulation results in Figure \ref{fig:sim_2}.
\begin{figure}
    \centering
    \includegraphics[scale=0.37]{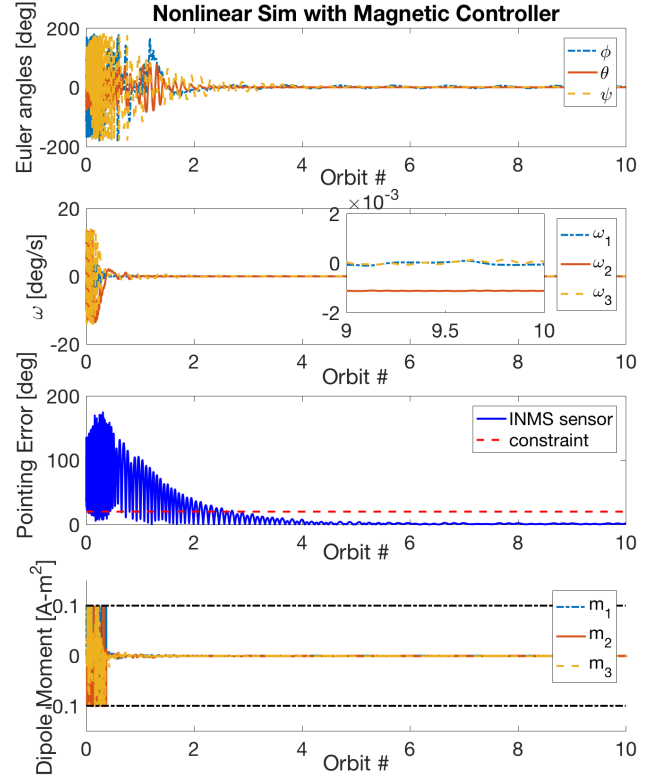}
    \caption{Simulation results using the predictive control design based on the relinearized dynamics, to include the panels, with the same constant disturbance torque as in Figure \ref{fig:dist}.}
    \label{fig:sim_mpc}
\end{figure}
The predictive controller provides faster convergence than the LQR controller, at the added cost of additional computation complexity and power consumption that may present challenges to a resource-limited cubesat platform.

\section{Conclusion}
\label{sec:concl}
The paper described the design of pointing controllers to enable the QB50 satellite's scientific mission. The two LQR controllers exploit the magnetic torque rod actuators to regulate both attitude and angular velocity states, and the second such controller complements the augmented passive drag panels. Both controllers have been shown to provide convergence to the desired pointing configuration in nonlinear model simulations; the second controller, however, has demonstrated greater robustness to unmodeled disturbance torques. A finite-horizon predictive controller is shown to provide faster convergence than the LQR controllers while maintaining the robustness to disturbance torques.

\section{Acknowledgment}
\label{sec:ack}
The authors would like to thank Professor Aaron Ridley and Mr. Patrick McNally for giving us an opportunity to contribute to QB50 development. We acknowledge them and the University of Michigan QB50 Team for numerous discussions on the mission and control system design requirements.

\section*{Appendix A}
\label{sec:appA}
\subsection*{Derivation of Linearized Equations of Motion}
\label{sec:d_leom}
We repeat (\ref{eqn:kinematics}) and (\ref{eqn:dynamics}), which express the full combined kinematics and dynamics in terms of the state variables:
\begin{equation}
    \begin{bmatrix}
        \dot{\phi} \\ \dot{\theta} \\ \dot{\psi}
    \end{bmatrix}
    = C^{-1}_{\phi \theta}
    \left(
    \begin{bmatrix}
        \omega_1 \\ \omega_2 \\ \omega_3
    \end{bmatrix}
    + n
    \begin{bmatrix}
        c_{\theta} s_{\psi} \\
        s_{\phi} s_{\theta} s_{\psi} + c_{\phi} c_{\psi} \\
        c_{\phi} s_{\theta} s_{\psi} - s_{\phi} c_{\psi}
    \end{bmatrix}
    \right),
\label{eqn:kinematics_a}
\end{equation}
\begin{equation}
    \mathbf{J}^{\mathcal{B}c}_b\dot{\boldsymbol{\omega}}^{bg}_b + S[\boldsymbol{\omega}^{bg}_b]\mathbf{J}^{\mathcal{B}c}_b\boldsymbol{\omega}^{bg}_b = \boldsymbol{\tau}^{\mathcal{B}c}_b.
\label{eqn:dynamics_a}
\end{equation}
where
\[
    C^{-1}_{\phi \theta} = \left( \frac{1}{c_\theta} \right)
    \begin{bmatrix}
        c_\theta & s_\phi c_\theta & c_\phi s_\theta \\
        0 & c_\phi c_\theta & -s_\phi c_\theta \\
        0 & s_\phi & c_\phi
    \end{bmatrix}.
\]

\subsubsection*{Linearized Kinematics}
\label{sec:lin_kin_a}
By design, the desired equilibrium state is such that [$\phi$, $\theta$, $\psi$] = [0, 0, 0], thus we can use small angle approximations $(c_{\phi} \approx 1$, $s_{\phi} \approx \phi)$ to simplify (\ref{eqn:kinematics_a}),
\begin{equation}
    \begin{bmatrix}
        \dot{\phi} \\ \dot{\theta} \\ \dot{\psi}
    \end{bmatrix}
    \approx
    \begin{bmatrix}
        1 & \phi \theta & \theta \\
        0 & 1 & -\phi \\
        0 & \phi & 1
    \end{bmatrix}
    \left(
    \begin{bmatrix}
        \omega_1 \\ \omega_2 \\ \omega_3
    \end{bmatrix}
    + n
    \begin{bmatrix}
        \psi \\ \phi \theta \psi + 1 \\ \theta \psi - \phi
    \end{bmatrix}
    \right).
\label{eqn:kinematics_a_sma1}
\end{equation}
Applying the equilibrium values for the Euler angles and Euler angle rates, we can easily verify that the equilibrium angular velocity values must satisfy [$\omega_1$, $\omega_2$, $\omega_3$] = [0, $-n$, 0]. Expanding (\ref{eqn:kinematics_a_sma1}):
\begin{equation}
    \begin{split}
        \dot{\phi} &= \omega_1 + \phi \theta \omega_2 + \theta \omega_3 + n \psi(1 + \phi^2 \theta^2 + \theta^2), \\
        \dot{\theta} &= \omega_2 - \phi \omega_3 + n(1 + \phi^2), \\
        \dot{\psi} &= \phi \omega_2 + \omega_3 + n(1 + \phi^2) \theta \psi.
    \end{split}
\label{eqn:kinematics_a_sma2}
\end{equation}
Then, taking the partial derivative in each state variable, we get the following sets of equations:
\begin{equation}
\begin{split}
    \frac{\partial \dot{\phi}}{\partial \phi} &= \theta \omega_2 + 2n \phi \theta^2 \psi, \\ \frac{\partial \dot{\phi}}{\partial \theta} &= \phi \omega_2 + 2n(1 + \phi^2)\theta \psi + \omega_3, \\ \frac{\partial \dot{\phi}}{\partial \psi} &= n(1 + \phi^2 \theta^2 + \theta^2), \\
    \frac{\partial \dot{\phi}}{\partial \omega_1} &= 1, \ \ \frac{\partial \dot{\phi}}{\partial \omega_2} = \phi \theta, \ \ \frac{\partial \dot{\phi}}{\partial \omega_3} = \theta,
\end{split}
\end{equation}
\begin{equation}
\begin{split}
    \frac{\partial \dot{\theta}}{\partial \phi} &= -\omega_3 + 2n\phi, \\ \frac{\partial \dot{\theta}}{\partial \theta} &= \frac{\partial \dot{\theta}}{\partial \psi} = \frac{\partial \dot{\theta}}{\partial \omega_1} = 0, \\ \frac{\partial \dot{\theta}}{\partial \omega_2} &= 1, \ \ \frac{\partial \dot{\theta}}{\partial \omega_3} = -\phi,
\end{split}
\end{equation}
\begin{equation}
\begin{split}
    \frac{\partial \dot{\psi}}{\partial \phi} &= \omega_2 + 2n\phi \theta \psi, \\ \frac{\partial \dot{\psi}}{\partial \theta} &= n(1 + \phi^2) \psi, \\ \frac{\partial \dot{\psi}}{\partial \psi} &= n(1+\phi^2) \theta, \\
    \frac{\partial \dot{\psi}}{\partial \omega_1} &= 0, \ \ \frac{\partial \dot{\psi}}{\partial \omega_2} = \phi, \ \ \frac{\partial \dot{\psi}}{\partial \omega_3} = 1.
\end{split}
\end{equation}

We complete the linearization of the kinematics by plugging the equilibrium values into the partial derivatives of the system:
\begin{equation}
    \begin{bmatrix}
        \dot{\phi} \\ \dot{\theta} \\ \dot{\psi}
    \end{bmatrix}
    =
    \begin{bmatrix}
        0 & 0 & n & 1 & 0 & 0 \\
        0 & 0 & 0 & 0 & 1 & 0 \\
        -n & 0 & 0 & 0 & 0 & 1
    \end{bmatrix}
    \begin{bmatrix}
        \phi \\ \theta \\ \psi \\ \omega_1 \\ \omega_2 + n \\ \omega_3
    \end{bmatrix}.
\label{eqn:kinematics_a_lin}
\end{equation}

\subsubsection*{Linearized Dynamics}
\label{sec:lin_dyn_a}
We are treating the case of the ideal, uncontrolled dynamics, so the only external torque effect to consider is the gravity gradient. Thus, we can replace the $\boldsymbol{\tau}_b^{\mathcal{B}c}$ term above with the equivalent $\boldsymbol{\tau}_b^{gg}$. Then, the dynamics can be be expressed as:
\begin{equation}
    \boldsymbol{J \dot{\omega}} + \text{S}\left[ \boldsymbol{\omega} \right]
    \boldsymbol{J \omega} = 3 n^2
    \begin{bmatrix}
        -(J_2 - J_3) c_\phi s_\phi c_\theta^2 \\
        (J_3 - J_1) c_\phi c_\theta s_\theta \\
        (J_1 - J_2) s_\phi c_\theta s_\theta
    \end{bmatrix}.
\label{eqn:dyn_a}
\end{equation}
This equation simplifies to:
\begin{equation}
    \begin{split}
        \dot{\omega}_1 = J_{23} ( \omega_2 \omega_3 - 3 n^2 c_\phi s_\phi c_\theta^2 ), \\
        \dot{\omega}_2 = J_{31} ( \omega_3 \omega_1 + 3 n^2 c_\phi c_\theta s_\theta ), \\
        \dot{\omega}_3 = J_{12} ( \omega_1 \omega_2 + 3 n^2 s_\phi c_\theta s_\theta ),
    \end{split}
\end{equation}
where $J_{12} := (J_1 - J_2)/J_3$, $J_{31} := (J_3 - J_1)/J_2$, and $J_{23} := (J_2 - J_3)/J_1$. \\\\
As with the kinematics, we now apply the small-angle approximations for the Euler angle terms:
\begin{equation}
    \begin{split}
        \dot{\omega}_1 &= J_{23} ( \omega_2 \omega_3 - 3 n^2 \phi ), \\
        \dot{\omega}_2 &= J_{31} ( \omega_3 \omega_1 + 3 n^2 \theta ), \\
        \dot{\omega}_3 &= J_{12} ( \omega_1 \omega_2 + 3 n^2 \phi \theta ).
    \end{split}
\end{equation}
We form the Jacobian of the small-angle system by taking the first-order partial derivatives in each state variable:
\begin{equation}
\begin{split}
    \frac{\partial \dot{\omega}_1}{\partial \phi} &= -3 n^2 J_{23}, \\ \frac{\partial \dot{\omega}_1}{\partial \theta} &= \frac{\partial \dot{\omega}_1}{\partial \psi} = \frac{\partial \dot{\omega}_1}{\partial \omega_1} = 0, \\ \frac{\partial \dot{\omega}_1}{\partial \omega_2} &= J_{23} \omega_3, \ \ \frac{\partial \dot{\omega}_1}{\partial \omega_3} = J_{23} \omega_2.
\end{split}
\end{equation}
\begin{equation}
\begin{split}
    \frac{\partial \dot{\omega}_2}{\partial \phi} &= \frac{\partial \dot{\omega}_2}{\partial \psi} = \frac{\partial \dot{\omega}_2}{\partial \omega_2} = 0, \\ \frac{\partial \dot{\omega}_2}{\partial \theta} &= 3 n^2 J_{31}, \\ \frac{\partial \dot{\omega}_2}{\partial \omega_1} &= J_{31} \omega_3, \ \ \frac{\partial \dot{\omega}_2}{\partial \omega_3} = J_{31} \omega_1,
\end{split}
\end{equation}
\begin{equation}
\begin{split}
    \frac{\partial \dot{\omega}_3}{\partial \phi} &= 3 n^2 J_{12} \theta, \ \ \frac{\partial \dot{\omega}_3}{\partial \theta} = 3 n^2 J_{12} \phi, \\ \frac{\partial \dot{\omega}_3}{\partial \psi} &= \frac{\partial \dot{\omega}_3}{\partial \omega_3} = 0, \\ \frac{\partial \dot{\omega}_3}{\partial \omega_1} &= J_{12} \omega_2, \ \ \frac{\partial \dot{\omega}_3}{\partial \omega_2} = J_{12} \omega_1,
\end{split}
\end{equation}
We complete the linearization by substituting the equilibrium values into the partial derivatives:
\begin{equation}
    \begingroup 
    \setlength\arraycolsep{1.3pt}
    \begin{bmatrix}
        \dot{\omega}_1 \\ \dot{\omega}_2 \\ \dot{\omega}_3
    \end{bmatrix}
    =
    \begin{bmatrix}
        -3 n^2 J_{23} & 0 & 0 & 0 & 0 & -n J_{23} \\
        0 & 3 n^2 J_{31} & 0 & 0 & 0 & 0 \\
        0 & 0 & 0 & -n J_{12} & 0 & 0
    \end{bmatrix}
    \begin{bmatrix}
        \phi \\ \theta \\ \psi \\ \omega_1 \\ \omega_2 + n \\ \omega_3
    \end{bmatrix}
    \endgroup
\label{eqn:dynamics_a_lin}
\end{equation}
The combined linearized kinematics and dynamics equations can now be expressed by the $A_c$ matrix that appears in (\ref{eqn:full_Ac}).

\section*{Appendix B}
\label{sec:appB}
\subsection{Algebraic Riccati Equation Solution Algorithm Formulation}
\label{sec:alg_form}
For implementation of the LQR controller in the satellite, we chose to follow a zero-order hold discrete-time formulation identical to (\ref{eqn:zoh}), with equations repeated here for convenience. Let $t \in \mathbb{R}_{\geq0}$ be the current time instant and $A_c$, $B_c(t)$ be the continuous-time dynamics defined in (\ref{eqn:full_Ac}), (\ref{eqn:full_Bc}). For $\Delta t > 0$, the discrete-time model predicts the state $\mathbf{x}_{k}$ at time $t+k\Delta t, k \in \mathbb{Z}_{\geq 0}$, according to the following model with the ``frozen-in-time'' magnetic field:
\begin{equation}
\begin{split}
    \mathbf{x}_{k+1} & = A_d \mathbf{x}_{k} + B_d(t) \mathbf{u}_{k}, \\
    A_d & = e^{A_c \Delta t}, \\
    B_d(t) & = - A_c^{-1} \left( \mathbf{I}_6 - A_d \right) B_c(t), \\
    \mathbf{x}_0 &= \mathbf{x}(t).
\end{split}
\label{eqn:zohB}
\end{equation}
The pair $(A_d, B_d(t))$ can be verified to be controllable for all $t$ for our orbit and choices of $\Delta t$. For the difference equation (\ref{eqn:zoh}), we define the infinite-horizon cost functional $J$:
\begin{equation}
    J = \sum_{k = 0}^{\infty} \left( \mathbf{x}_k^T Q \mathbf{x}_k + \mathbf{u}_k^T R \mathbf{u}_k \right),
\end{equation}
where $R=R^T \in \mathbb{R}^3$ is a positive definite matrix, and $Q=Q^T \in \mathbb{R}^6$ is a positive semi-definite matrix satisfying the usual detectability assumption. Then, the optimal feedback control sequence $\mathbf{u}_k = -K(t) \mathbf{x}_k$ that minimizes $J$ has the solution:
\begin{equation}
    K(t) = \left( R + B_d(t)^T P(t) B_d(t) \right)^{-1} B_d(t)^T P(t) A_d,
\label{eqn:kdB}
\end{equation}
\begin{equation}
    P(t) = A_d^T P(t) A_d + Q - A_d^T P(t) B_d(t) K(t).
\label{eqn:dareB}
\end{equation}
Note that (\ref{eqn:dare}) can have multiple solutions; the $P(t)$ of interest to us is the unique positive definite solution. Also note that $B_d(t)$ changes throughout the orbit, thus (\ref{eqn:dare}) is solved at different instants $t$ in time and the gain $K(t)$ in (\ref{eqn:kd}) is time-varying. The control $\mathbf{u}(t+\sigma) = K(t) \mathbf{x}(t)$ is applied for $0 \leq \sigma < \Delta t$ and then recomputed.
\\
\subsubsection{Gain Computation}
Upon generating the discrete-time model in (\ref{eqn:zohB}), we define the following two matrices:
\begin{equation}
\begin{split}
    N &:=
    \begin{bmatrix}
        A_d & \mathbf{0}_{6 \times 6} \\
        -Q & \mathbf{I}_6
    \end{bmatrix}, \\
    L &:=
    \begin{bmatrix}
        \mathbf{I}_6 & B_d R^{-1} B_d^T \\
        \mathbf{0}_{6 \times 6} & A_d^T
    \end{bmatrix},
\end{split}
\end{equation}
where $A_d \in \mathbb{R}^{6 \times 6}$ and $B_d = B_d(t) \in \mathbb{R}^{6 \times 3}$. From these, we form the Hamiltonian matrix
\begin{equation}
    H = \left( N + L \right) ^{-1} \left( N - L \right).
\end{equation}
We require the positive square root of $H^2$, and employ a Newton-Raphson iteration process to compute it. We begin with an initial guess of $S_0 = \mathbf{I}_6$, the identity matrix, and then iterate according to the following scheme:
\begin{equation}
    S_{k+1} = 0.5 \left( S_k + S_k^{-1} H^2 \right).
\end{equation}
Upon reaching desired convergence in $S_k$, we can then extract the unique positive-definite solution $P(t)$ to (\ref{eqn:dare}) from the first column of the matrix
\begin{equation}
    \begin{bmatrix}
        X_1 & \sim \\ X_2 & \sim
    \end{bmatrix}
    = H - S,
\end{equation}
as
\begin{equation}
    P(t) = X_2 X_1^{-1}.
\end{equation}
Given $P(t)$, we then compute
\begin{equation}
    K(t) = \left( R + B_d(t)^T P(t) B_d(t) \right)^{-1} B_d(t)^T P(t) A_d.
\end{equation}
This solution algorithm is suitable for use even on our resource-limited cubesat platform, as the necessary matrix inverses can be computed very efficiently through a decomposition scheme; in our case, we applied an LU-decomposition, as described in Anton \& Rorres\cite{LinAlg}. Furthermore, since the solution between time steps does not change significantly in this slowly-varying system, after the first gains are computed then future solutions can be ``warm started'', i.e., the initialization $S_0$ takes the value of the solution from the previous time step. This has the benefit of reducing the number of iterations required to obtain convergence. Once the solution matrix $P(t)$ is found, it is then a simple matter to compute the LQR gain matrix $K(t)$. Note that $A_d$, $Q$, and $R^{-1}$ are invariant and do not depend on the magnetic field, hence are precomputed and stored in memory.

\begin{figure}
    \centering
    \includegraphics[scale = 0.19]{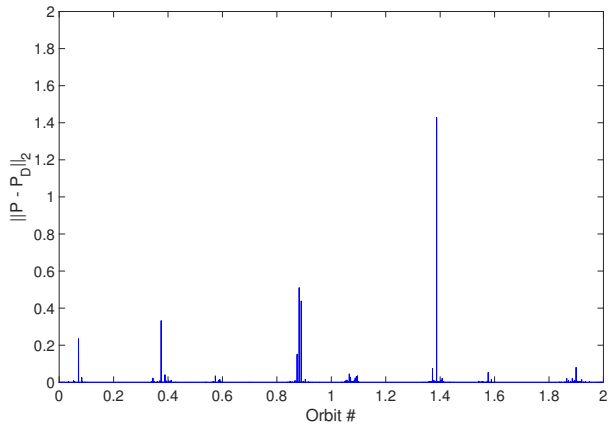}
    \caption{2-norm of the difference between the exact solution to the Discrete-Time Algebraic Riccati Equation and the approximate solution returned by the algorithm.}
    \label{fig:PmPD}
\end{figure}

We measure the accuracy of the algorithm's output by taking the 2-norm of the difference between the approximate solutions, denoted $P_D(t)$, with the exact solutions $P(t)$ returned by Matlab's \textit{dlqr} command over two simulated orbits. See Figure \ref{fig:PmPD}.
The norms are generally close to zero, indicating very strong agreement between the approximate and exact solutions, validating the use of the described algorithm.

\section*{Appendix C}
\label{sec:appC}
\subsection*{Proof of Theorem 1}
LTV controllability holds if the controllability matrix analogue $\mathcal{K} = [ K_0, \ K_1, \ K_2 ]$ has full row rank for some $t_c \in [0, \ T]$\cite{Yang}, where $K_j = \frac{\partial^j}{\partial \tau^j} [ \Phi(t, \tau) B(\tau) ]_{\tau = t}$ and $\Phi(t, \tau)$ is the state transition matrix for the $A$ matrix. We form this matrix and show that non-singularity holds for $t_c = T/4$ under the above conditions.

Suppose the satellite's orbit aligns with the magnetic Equator, so that $i_m = 0$. Then, for all $t$, $b_1(t) = b_3(t) = 0$ and $b_2(t)$ is constant. Then, $B$ is constant and takes the form,
\begin{equation}
    B = \begin{bmatrix} & \mathbf{0}_{3 \times 3} & \\ b_{\alpha} & 0 & 0 \\ 0 & 0 & 0 \\ 0 & 0 & b_{\gamma} \end{bmatrix}.
\end{equation}

In this case, the lack of controllability can be shown by forming the standard controllability matrix for LTI systems, $\mathcal{C} = \begin{bmatrix} B & AB & A^2 B \end{bmatrix} \in \mathbb{R}^{6 \times 9}$, and showing that it does not have full row rank. Each of $B$, $AB$, and $A^2 B$ has all zero entries in its fifth row, thus $\mathcal{C}$ has all zero entries in its fifth row. Thus, the row rank of $\mathcal{C}$ is at most $5 < 6$, $\mathcal{C}$ does not have full row rank, and the associated LTI system is not controllable.

Now assume that the satellite's orbit is inclined relative to the magnetic Equator. $B$ is time-varying, thus we must form the time-varying controllability analogue,
\begin{equation}
    \mathcal{K} = 
    \begingroup 
    \setlength\arraycolsep{3pt} 
    \begin{bmatrix} B(t) & \dot{B}(t) - A B(t) & A^2 B(t) - 2 A \dot{B}(t) + \ddot{B}(t) \end{bmatrix} \endgroup ,
\end{equation}
and show that the matrix $\mathcal{K}$ has full row rank for some $t_c \in [0, \ T]$. For convenience, select $t_c = T/4$, as then $nt_c = \pi/2$ and the trigonometric terms in the magnetic field approximation simplify greatly. Following Yang \cite{Yang}, we express the $A$ and $B$ matrices as follows,
\begin{equation}
    A = \begin{bmatrix} \Sigma_1 & \mathbf{I}_3 \\ \Lambda & \Sigma_2 \end{bmatrix}, \ B = \begin{bmatrix} \boldsymbol{0}_{3 \times 3} \\ B_2 \end{bmatrix},
\end{equation}
so that $\mathcal{K}$ can then be expressed as
\begin{equation}
    \mathcal{K} = 
    \begingroup 
    \setlength\arraycolsep{2.5pt} 
    \begin{bmatrix} \boldsymbol{0}_{3 \times 3} & -B_2 & \left(\Sigma_1 + \Sigma_2\right)B_2 - 2\dot{B}_2 \\ B_2 & -\Sigma_2 B_2 + \dot{B}_2 & \left(\Lambda + \Sigma^2_2\right) B_2 - 2 \Sigma_2 \dot{B}_2 + \ddot{B}_2 \end{bmatrix}
    \endgroup .
\end{equation}

We now look for a submatrix of $\mathcal{K}$ in $\mathbb{R}^{6 \times 6}$ that is non-singular. We can simplify $\mathcal{K}$ with a row reduction by premultiplying the top row by $-\Sigma_2$ and adding the result to the second row to get
\begin{equation}
\begin{split}
    \mathcal{K}_2 &= \begin{bmatrix} \boldsymbol{0}_{3 \times 3} & -B_2 & \left(\Sigma_1 + \Sigma_2\right)B_2 - 2\dot{B}_2 \\ B_2 & \dot{B}_2 & \left(\Lambda - \Sigma_2 \Sigma_1\right) B_2 + \ddot{B}_2 \end{bmatrix} \\ &= \begin{bmatrix} \boldsymbol{0}_{3 \times 3} & -B_2 & M_1 \\ B_2 & \dot{B}_2 & M_2 \end{bmatrix} .
\end{split}
\end{equation}

$\mathcal{K}_2$ having a non-singular submatrix is equivalent to $\mathcal{K}$ having one, so we now work with the new matrix instead. At the chosen time instant, we have
\begin{equation}
    \begin{split}
        b_1(t_c) &= 0, \ \dot{b}_1(t_c) = -n \frac{\mu_f}{a^3} s_{i_m}, \ \ddot{b}_1(t_c) = -n^2 \frac{\mu_f}{a^3} s_{i_m} \\
        b_2(t_c) &= -\frac{\mu_f}{a^3} c_{i_m}, \ \dot{b}_2(t_c) = 0, \ \ddot{b}_2(t_c) = 0 \\
        b_3(t_c) &= 2 \frac{\mu_f}{a^3} s_{i_m}, \ \dot{b}_3(t_c) = 0, \ \ddot{b}_3(t_c) = -n^2 \frac{\mu_f}{a^3} s_{i_m}.
    \end{split}
\end{equation}

Define $p_1 = (\mu_f/a^3) s_{i_m}$, $p_2 = (\mu_f/a^3) c_{i_m}$, $p_3 = 2 (\mu_f/a^3) s_{i_m} = 2 p_1$, and $p_{jk} = p_j / J_k$. Then, $\mathcal{K}_2$ takes the form,
\begin{equation}
    \mathcal{K}_2 = 
    \begingroup 
    \setlength\arraycolsep{1.5pt}
    \begin{bmatrix} 
        0 & 0 & 0 & 0 & -p_{31} & p_{21} & m_{11} & 0 & 0 \\
        0 & 0 & 0 & p_{32} & 0 & 0 & 0 & 0 & m_{12} \\
        0 & 0 & 0 & -p_{23} & 0 & 0 & 0 & m_{13} & m_{14} \\
        0 & p_{31} & -p_{21} & 0 & 0 & 0 & 0 & m_{21} & m_{22} \\
        -p_{32} & 0 & 0 & 0 & 0 & n p_{12} & m_{23} & 0 & 0 \\
        p_{23} & 0 & 0 & 0 & -n p_{13} & 0 & m_{24} & 0 & 0
    \end{bmatrix}
    \endgroup,
\end{equation}
and we form a submatrix $\mathcal{K}_3$ from columns 1, 2, 4, 5, 7, and 8.
    \begin{equation}
        \mathcal{K}_3 = 
        \begin{bmatrix}
            0 & 0 & 0 & -p_{31} & m_{11} & 0 \\
            0 & 0 & p_{32} & 0 & 0 & 0 \\
            0 & 0 & -p_{23} & 0 & 0 & m_{13} \\
            0 & p_{31} & 0 & 0 & 0 & m_{21} \\
            -p_{32} & 0 & 0 & 0 & m_{23} & 0 \\
            p_{23} & 0 & 0 & -n p_{13} & m_{24} & 0
        \end{bmatrix},
    \end{equation}
which we can row-reduce to,
\begin{equation}
    \mathcal{K}_4 = 
    \begin{bmatrix}
        0 & 0 & 0 & -p_{31} & m_{11} & 0 \\
        0 & 0 & p_{32} & 0 & 0 & 0 \\
        0 & 0 & 0 & 0 & 0 & m_{13} \\
        0 & p_{31} & 0 & 0 & 0 & 0 \\
        -p_{32} & 0 & 0 & 0 & m_{23} & 0 \\
        p_{23} & 0 & 0 & -n p_{13} & m_{24} & 0
    \end{bmatrix},
\end{equation}
which has determinant det$(\mathcal{K}_4) = (-p_{32})(-m_{13})(p_{31})(m_{24}p_{31}p_{32} - m_{23}(-p_{31})p_{23} + m_{11}p_{32} nb_{13})$. The LTV system is then controllable if this determinant is nonzero. Each $p$ term is already nonzero, thus we seek $m_{13} \neq 0$ and det$(\mathcal{K}_5) = m_{24}p_{31}p_{32} + m_{23}p_{31}p_{23} - m_{11}p_{32} np_{13} \neq 0$. We have:
\begin{equation}
    \begin{split}
        m_{13} &= \frac{-2n p_1}{J_3} - \frac{n(-J_3 + J_2 - J_1) p_3}{J_1 J_3} \\
               &= \left( \frac{-2n p_1}{J_1 J_3} \right) (J_1 - J_3 + J_2 - J_1) \\
               &= \left( \frac{-2n p_1}{J_1 J_3} \right) (J_2 - J_3) \\
               &\neq 0
    \end{split}
\end{equation}
which, as $n \neq 0$ and $p_1 \neq 0$, implies that $J_2 \neq J_3$.

Next,
\begin{equation}
    \begin{split}
        m_{11}p_{32} np_{13} &= \frac{-n p_1}{J_3} \frac{p_3}{J_2} \frac{p_2}{J_3} \left(\frac{nJ_1 + n(J_3 - J_2)}{J_1}\right) \\
        &= \left(\frac{p_2 p_3^2 n^2}{J_1 J_2 J_3}\right) \left(\frac{-J_1 + J_2 - J_3}{2 J_3} \right) \\
        m_{23}p_{31}p_{23} &= \frac{p_3}{J_2} \frac{p_3}{J_1} \frac{p_2}{J_3} \left( 3n^2 \frac{J_3 - J_1}{J_2} + n^2 \frac{\Gamma}{J_2} - n^2 \right) \\
        &= \left(\frac{p_2 p_3^2 n^2}{J_1 J_2 J_3}\right) \left(\frac{3(J_3-J_1) + \Gamma - J_2}{J_2}\right) \\
        m_{24}p_{31}p_{32} &= \frac{-p_2}{J_3} \frac{p_3}{J_1} \frac{p_3}{J_2} \left( n^2 \frac{\Gamma}{J_3} - n^2 \frac{J_2 - J_1}{J_3} \right) \\
        &= \left(\frac{p_2 p_3^2 n^2}{J_1 J_2 J_3}\right) \left(\frac{J_2 - J_1 - \Gamma}{J_3}\right),
    \end{split}
\end{equation}
thus, the condition
\begin{equation}
    \text{det} (\mathcal{K}_5) \neq 0
\end{equation}
implies, after some regrouping, that
\begin{equation}
    J_3 \left( 6(J_3 - J_1) + 2 \Gamma \right) \neq J_2 \left(J_1 - J_2 + J_3 - 2 \Gamma \right).
\end{equation}

\bibliographystyle{unsrt}  
\bibliography{QB50}

\end{document}